\documentclass[12pt]{amsart}
\numberwithin{equation}{section}

\usepackage{amssymb}
\def\cb{{\mathcal B}}

\def\ch{{\mathcal H}}

\def\ck{{\mathcal K}}

\def\cs{{\mathcal S}}

\def\ga{{\mathfrak A}} 
\def\gb{{\mathfrak B}}

\def\gam{{\mathfrak M}}
\def\gn{{\mathfrak N}}

\def\bc{{\mathbb C}}
\def\bbf{{\mathbb F}}

\def\bn{{\mathbb N}}

\def\bt{{\mathbb T}}

\def\bz{{\mathbb Z}}

\def\a{\alpha}
\def\b{\beta}
\def\g{\gamma} \def\G{\Gamma}
  \def\D{\Delta}

\def\eps{\varepsilon}

\def\r{\rho}
\def\s{\sigma} \def\S{\Sigma}
\def\t{\tau}
\def\f{\varphi} \def\F{\Phi}
  
\def\om{\omega} \def\Om{\Omega}

\def\id{\hbox{id}}

\newtheorem{thm}{Theorem}[section]
\newtheorem{lem}[thm]{Lemma}

\newtheorem{cor}[thm]{Corollary}
\newtheorem{prop}[thm]{Proposition}
\newtheorem{defin}[thm]{Definition}

\def\ad{\mathop{\rm Ad}}

\newcommand{\nn}{\nonumber}

\newcommand{\ots}[1]{\otimes_{\mathop{\rm #1}}}
\def\oots{\overline{\otimes}}

\begin{document}

\title[entangled ergodic theorem, diagonal measures]
{the entangled ergodic theorem and an ergodic theorem for quantum ``diagonal measures''}
\author{Francesco Fidaleo}
\address{Francesco Fidaleo,
Dipartimento di Matematica,
Universit\`{a} di Roma Tor Vergata, 
Via della Ricerca Scientifica 1, Roma 00133, Italy} \email{{\tt
fidaleo@mat.uniroma2.it}}


\begin{abstract}
Let $U$ be a unitary operator acting on the Hilbert space $\ch$, 
$\a:\{1,\dots, 2k\}\mapsto\{1,\dots, k\}$ a pair--partition, and 
finally $A_{1},\dots,A_{2k-1}\in\cb(\ch)$.
We show that the ergodic average 
$$
\frac{1}{N^{k}}\sum_{n_{1},\dots,n_{k}=0}^{N-1}
U^{n_{\a(1)}}A_{1}U^{n_{\a(2)}}\cdots 
U^{n_{\a(2k-1)}}A_{2k-1}U^{n_{\a(2k)}}
$$
converges in the strong 
operator topology when $\ch$ is generated by the eigenvectors of $U$, 
that is when the dynamics induced by the unitary $U$ on 
$\ch$ is almost periodic. This result improves the known ones 
relative to the 
entangled ergodic theorem. 
We also prove the noncommutative version of the ergodic result 
of H. Furstenberg relative to diagonal measures. This implies that
${\displaystyle
\frac{1}{N}\sum_{n=0}^{N-1}
U^{n}AU^{n}}$ converges in the strong operator topology for other 
interesting situations where 
the involved unitary operator does not generate an almost periodic 
dynamics, and the operator $A$ is noncompact.

\end{abstract}

\maketitle

\section{introduction}

The investigation of ergodic properties of classical dynamical systems 
has a long history. As an example, we mention the well--known ergodic hypothesis 
(cf. e.g. \cite{LL}, Section 4) which can be viewed as a justification of the 
microcanonical 
distribution in statistical mechanics. 
We refer the reader to \cite{RS} for a nice introduction, and
the monograph \cite{AA} for the basic 
results and further details.

Recently, the ergodic theory of noncommutative dynamical systems has 
been an impetuos growth in relation to the natural applications to 
quantum (statistical) physics. In view to other potential 
applications, it is of interest to understand among the various 
ergodic properties, which ones survive by passing from the classical 
to the quantum case. We mention the pivotal paper \cite{NSZ}, where such a program is 
carried out for some basic recurrence, as well as multiple mixing properties.

Notice that 
it is in general unclear what should be the right quantum counterpart 
of a classical ergodic property. For example, the reader can compare   
the property of the 
convergence to the equilibrium (i.e. ergodicity for an invariant 
state $\om$)
$$
\lim_{N\to+\infty}\frac{1}{N}\sum_{n=0}^{N-1}\om\left(B^{*}\a^{n}(A)B\right)
=\om\left(B^{*}B\right)\om(A)
$$
suggested by the quantum physics, with the standard notion of 
ergodicity
$$
\lim_{N\to+\infty}\frac{1}{N}\sum_{n=0}^{N-1}\om(A\a^{n}(B))
=\om(A)\om(B)\,.
$$
See \cite{FL}, Proposition 1.1 for further details.

A notion which is meaningful in quantum setting is that of 
entangled ergodic theorem, formulated in \cite{AHO} in 
connection with the central limit theorem for suitable sequences of elements of the group 
$C^{*}$--algebra of the free group $\bbf_{\infty}$ on infinitely many generators. 

The entangled ergodic theorem was clearly formulated in \cite{L}.
Namely, let $U$ be 
a unitary operator acting on the Hilbert space $\ch$, and for $m\geq k$, 
$\a:\{1,\dots,m\}\mapsto\{1,\dots, k\}$ a partition of the 
set $\{1,\dots,m\}$ in $k$ parts. 
The entangled ergodic theorem concerns 
the convergence in the strong, or merely weak operator topology, 
of the multiple Cesaro mean
\begin{equation}
\label{0}
\frac{1}{N^{k}}\sum_{n_{1},\dots,n_{k}=0}^{N-1}
U^{n_{\a(1)}}A_{1}U^{n_{\a(2)}}\cdots 
U^{n_{\a(m-1)}}A_{m-1}U^{n_{\a(m)}}\,,
\end{equation}
$A_{1},\dots,A_{m-1}$ being bounded operators acting on $\ch$.

Expressions like \eqref{0} naturally appear also in \cite{NSZ}
relatively to the study of the 
multiple mixing. Namely, suppose that the dynamics of a (concrete) 
dynamical system is unitarily implemented by the unitary $U$, and the 
vector $\Om$ is invariant under $U$.\footnote{Notice that this is 
always the case by considering the GNS covariant representation.} 
Then firstly in \cite{Fu}, and 
more recently in \cite{NSZ}, the behavior of the multiple correlations
\begin{equation}
\label{00}
\frac{1}{N}\sum_{n=0}^{N-1}\om\left(A_{0}\ad{}\!_{U}^{n}(A_{1})
\ad{}\!_{U}^{2n}(A_{2})\right)
\equiv\frac{1}{N}\sum_{n=0}^{N-1}\left\langle U^{n}A_{1}U^{n}A_{2}
\Om,A^{*}_{0}\Om\right\rangle
\end{equation}
has been studied in connection with the $(1,2)$--multiple mixing or merely ergodicity. 
Notice that \eqref{00} is the particular case of \eqref{0} relative to the (trivial) 
pair--partition of two elements. Just by considering the simplest case 
of the  
partition of the empty set, the limit of the 
Cesaro mean in \eqref{0} reduces itself to the well--known mean ergodic theorem 
due to John von Neumann (cf. \cite{RS}) 
\begin{equation}
\label{jvn}
\mathop{\rm s\!-\!lim}_{N\to+\infty}\frac{1}{N}\sum_{n=0}^{N-1}U^{n}=E_{1}\,,
\end{equation}
$E_{1}$ being the selfadjoint projection onto the eigenspace of the 
invariant vectors for $U$.

Some applications of the entangled ergodic theorem are discussed in 
\cite{F1}. Apart from the other potential applications to the study of the ergodic 
properties of quantum dynamical systems, the entangled ergodic theorem 
is a fascinating self--contained mathematical problem. It is certainly true if the 
spectrum $\s(U)$ of $U$ is finite. Some very special cases for which 
it holds true are listed in \cite{L}. It was shown in \cite{F} that
the entangled ergodic theorem holds true in a sufficiently general 
situation, that is when the operators 
$A_{1},\dots,A_{m-1}$ in \eqref{0} are compact.

The first part of the present paper is devoted to prove 
the entangled ergodic theorem in the case when the unitary $U$ is 
almost periodic (i.e. when $\ch$ is generated by the eigenvectors of 
$U$) and $\a:\{1,\dots,2k\}\mapsto\{1,\dots, k\}$ a pair--partition, 
without any condition on the operators $A_{1},\dots,A_{2k-1}$. 
This result improves those in Section 3 of \cite{F} relative to the almost periodic 
case, where only very special pair--partitions were considered.

The entangled ergodic theorem is not yet available in the full 
generality. Then it is natural to address the problem to find other 
nontrivial cases for which it holds true.

Another situation of interest arises from the generalization to the 
noncommutative setting, of the ergodic theorem of H. Furstenberg relative 
to diagonal measure (cf. \cite{Fu, Fu1}). This is precisely the 
argument of the second part of the present paper. Namely, we prove an 
ergodic theorem relative to possibly noninvariant and nonnormal states, 
which is the generalization of Theorem 3.1 of \cite{Fu1} relative to 
the Abelian case. This allows us to prove the following result. Let 
$M$ be a von Neumann algebra equipped with the adjoint action of an ergodic unitary $U$,
and a standard vector $\Om$ which is invariant under $U$. Let $M'$ be the 
commutant von Neumann algebra of $M$. The state defined as
$$
A\otimes B\in M\otimes M'\mapsto\langle AB\Om,\Om\rangle\,,
$$
is precisely the quantum counterpart of the ``diagonal measure'' 
associated to the product state
$$
A\otimes B\in M\otimes M'\mapsto\langle A\Om,\Om\rangle
\langle B\Om,\Om\rangle\,.
$$
We show that the Cesaro mean
$$
\frac{1}{N}\sum_{n=0}^{N-1}U^{n}AU^{n}
$$
converges in the strong operator topology for each $A\in M\bigcup M'$.

\section{terminology, notations and basic results}

Let $X$, $Y$ be linear spaces. Their algebraic tensor product is 
denoted by $X\otimes Y$. If $\ch$ and $\ck$ are Hilbert spaces, the 
Hilbertian tensor product, that is the completion of $\ch\otimes \ck$
under the norm induced by the inner product 
$$
\langle x\otimes\xi,y\otimes\eta\rangle:=
\langle x,y\rangle\langle\xi,\eta\rangle\,,
$$
is denoted as $\ch\oots \ck$.

Let $\{A_{\a}\}_{\a\in J}\subset\cb(\ch)$ be a net consisting of 
bounded operators acting on the Hilbert space $\ch$. If it converges 
to $A\in\cb(\ch)$ in the weak operator topology, respectively strong
operator topology, we write respectively
$$
\mathop{\rm w\!-\!lim}_{\a}A_{\a}=A\,,\quad
\mathop{\rm s\!-\!lim}_{\a}A_{\a}=A\,.
$$

Let $U$ be a unitary operator acting on $\ch$. Consider the resolution of the identity  
$\{E(\D)\,:\,\D\,\,\text{Borel subset of}\,\,\bt\}$ of $U$ (cf. 
\cite{TL}, Section VII.7). Denote with an abuse of notation, 
$E_{z}:=E(\{z\})$. Namely, $E_{z}$ is nothing but the selfadjoint 
projection on the eigenspace corresponding to the eigenvalue 
$z$ in the unit circle $\bt$. 

The unitary $U$ is said to be {\it ergodic} if the 
fixed--point subspace $E_{1}\ch$ is one 
dimensional. By the mean ergodic theorem \eqref{jvn}, it is equivalent to the 
existence of a unit vector $\xi_{0}\in\ch$ such that 
$$
\lim_{N\to+\infty}\frac{1}{N}\sum_{n=0}^{N-1}U^{n}\xi=
\langle\xi,\xi_{0}\rangle\xi_{0}\,,
$$
or equivalently,
$$
\lim_{N\to+\infty}\frac{1}{N}\sum_{n=0}^{N-1}\langle U^{n}\xi,\eta\rangle=
\langle\xi,\xi_{0}\rangle\langle\xi_{0},\eta\rangle\,.
$$

The unitary $U$ is said to be {\it weakly mixing} if there exists 
a unit vector $\xi_{0}\in\ch$ such that 
$$
\lim_{N\to+\infty}\frac{1}{N}\sum_{n=0}^{N-1}\big|\langle U^{n}\xi,\eta\rangle
-\langle\xi,\xi_{0}\rangle\langle\xi_{0},\eta\rangle\big|=0\,.
$$
Of course, a weakly mixing unitary is ergodic. It is well--known that 
the {\it vice--versa} does not hold. Indeed, $U$ is ergodic if and 
only if $E_{1}\ch$ is one dimensional. It is weakly mixing if and only 
if in addition, $\s_{\mathop{\rm pp}}(U)=\{1\}$, $\s_{\mathop{\rm 
pp}}(U)$ being the pure point spectrum of $U$ (cf. \cite{RS}). See 
e.g. \cite{NSZ}. 

The unitary $U$ is said to be {\it almost periodic} if
$\ch=\ch_{\mathop{\rm ap}}^{U}$, $\ch_{\mathop{\rm ap}}^{U}$ being 
the closed subspace 
consisting of the vectors having relatively norm--compact orbit 
under $U$. It is seen in \cite{NSZ} that $U$ is almost periodic if and 
only if $\ch$ is generated by the eigenvectors of $U$.

Define
\begin{equation}
\label{asy}
\s_{\mathop{\rm pp}}^{\mathop{\rm a}}(U):=
\big\{z\in\s_{\mathop{\rm pp}}(U)\,:\,zw=1\,\text{for 
some}\,w\in\s_{\mathop{\rm pp}}(U)\big\}\,.
\end{equation}
It is immediate to verify that $\s_{\mathop{\rm pp}}^{\mathop{\rm a}}(U)$
is a subgroup of the unit circle $\bt$.

Let $\a:\{1,\dots,2k\}\mapsto\{1,\dots,k\}$ 
be a pair--partition of the set $\{1,\dots,k\}$. 
It is 
shown in Proposition 2.3 of \cite{F}, that the net
\begin{align*}
\bigg\{\sum_{z_{1},\dots,z_{k}\in F}&E_{z^{\#}_{\a(1)}}A_{1}E_{z^{\#}_{\a(2)}}
\cdots
E_{z^{\#}_{\a(2k-1)}}A_{2k-1}E_{z^{\#}_{\a(2k)}}\,\\
:\,& 
F\subset\s_{\mathop{\rm pp}}^{\mathop{\rm a}}(U)
\,\text{finite subsets}\,\bigg\}\subset\cb(\ch)
\end{align*}
converges in the weak operator topology to a bounded operator written symbolically as
\begin{align}
\label{symb}
&S_{\a;A_{1},\dots,A_{2k-1}}\\
=\sum_{z_{1},\dots,z_{k}\in\s_{\mathop{\rm pp}}^{\mathop{\rm a}}(U)}
&E_{z^{\#}_{\a(1)}}A_{1}E_{z^{\#}_{\a(2)}}\cdots
E_{z^{\#}_{\a(2k-1)}}A_{2k-1}E_{z^{\#}_{\a(2k)}}\nn\,.
\end{align}

More precisely, the pairs $z^{\#}_{\a(i)}$ are alternatively $z_{j}$ and
$\bar z_{j}$ whenever $\a(i)=j$, 
and finally
the sum is understood as the limit in the weak 
operator topology of the above mentioned net obtained by considering 
all the finite truncations of the r.h.s. of 
\eqref{symb}.\footnote{If for example, $\a$ is the pair--partition
$\{1,2,1,2\}$ of four elements, 
$$
S_{\a;A,B,C}=\sum_{z,w\in\s_{\mathop{\rm pp}}^{\mathop{\rm a}}(U)}
E_{z}AE_{w}BE_{\bar z}CE_{\bar w}\,.
$$}  

Fix a pair--partition $\b:\{1,\dots, 2k+2\}\mapsto\{1,\dots, k+1\}$. Let 
$k_{\b}\in\{1,\dots, 2k+1\}$ be the first element of the pair 
$\b^{-1}\big(\{k+1\}\big)$, and $\a_{\b}$ the pair--partition of
$\{1,\dots, 2k\}$ obtained 
by deleting $\b^{-1}\big(\{k+1\}\big)$ from $\{1,\dots, 2k+2\}$, and 
$k+1$ from $\{1,\dots, k+1\}$. 
Notice that, if $x\in\ch$ is an eigenvector of $U$ with eigenvalue 
$z_{0}$, then we obtain
\begin{equation}
\label{zaz}
S_{\b;A_{1},\dots,A_{2k+1}}x
=S_{\a_{\b};A_{1},\dots,A_{k_{\b}-1}
E_{\bar z_{0}}A_{k_{\b}+1},\dots,A_{2k}}A_{2k+1}x\,.
\end{equation}

For a (discrete) $C^*$--dynamical system we mean a triplet $\big(\ga,\a,\om\big)$ 
consisting of a $C^*$-algebra $\ga$, an automorphism $\a$ of $\ga$, 
and a state $\om\in\cs(\ga)$ invariant under the action of $\a$. 

A $C^*$--dynamical system $\big(\ga,\a,\om\big)$
is said to be {\it ergodic} if for each 
$A,B\in\ga$,\footnote{Notice that this definition of ergodicity for 
an invariant state differs 
from the standard one. Indeed, an invariant state $\f\in\cs(\ga)$ is 
said to be ergodic if it is extremal among all the states 
of $\ga$ which are invariant under the action of $\a$. It can be shown that if $\f$ is 
asymptotically Abelian (or for an even and graded--asymptotically 
Abelian state $\f$, when $\ga$ is a $\bz_{2}$--graded $C^*$-algebra), both definitions 
coincide. See e.g. \cite{F0}, Section 3 for further details.}
$$
\lim_{N\to+\infty}\frac{1}{N}\sum_{n=0}^{N-1}\om(A\a^{n}(B))=
\om(A)\om(B)\,.
$$
It is said to be {\it weakly mixing} if
$$
\lim_{N\to+\infty}\frac{1}{N}\sum_{n=0}^{N-1}\big|\om(A\a^{n}(B))-
\om(A)\om(B)\big|=0
$$
for each $A,B\in\ga$. 

Let $\big(\ch,\pi,U,\Om\big)$ be the GNS covariant representation 
(cf. \cite{T}, Section I.9)  
canonically associated to the dynamical system under 
consideration. Then $\big(\ga,\a,\om\big)$ is ergodic (respectively 
weakly mixing) if and only if $U$ is ergodic (respectively 
weakly mixing), see e.g. \cite{NSZ}. 

Let $s(\om)$ be the support of $\om$ in the bidual 
$\ga^{**}$. Then $s(\om)\in Z(\ga^{**})$ if and only if $\Om$ is 
separating for $\pi(\ga)''$, $Z(\ga^{**})$ being the centre of 
$\ga^{**}$ (see e.g. \cite{SZ}, Section 10.17).

Denote $M:=\pi(\ga)''$, and with an abuse of notation, $\a:=\ad_{U}$ 
the adjoint action of $U$ on $\cb(\ch)$. The commutant von Neumann 
algebra is $M'\equiv \pi(\ga)'$.  
For $z$ in $\bt$ denote
$$
M_z=\{A\in M: \a(A)=zA\}\,,\quad (M')_z=\{B\in M': \a(B)=zB\}\,.
$$

The following results are probably known to the experts. We provide 
their proof for the convenience of the reader.
\begin{prop}
\label{gnsc}
Let the $C^*$--dynamical system $\big(\ga,\a,\om\big)$ be such that 
$s(\om)$ is central. Then, with the previous notations, 
$$
\overline{M_z\Om}=\overline{(M')_z\Om}=E_{z}\ch\,,
$$
and we can choose an orthonormal basis 
$\{u^{z}_{\a_{z}}\}_{\a_{z}\in I_{z}}\subset M_z\Om$ (equivalently 
$\{v^{z}_{\b_{z}}\}_{\b_{z}\in J_{z}}\subset (M')_z\Om$) for $E_{z}\ch$.

In addition, $\s_{\mathop{\rm pp}}(U)=\s_{\mathop{\rm pp}}(U)^{-1}$, and if
$z\in\s_{\mathop{\rm pp}}(U)$,
$$
\bigg\{\frac{\big(A^{z}_{\a_{z}}\big)^{*}\Om}{\sqrt{\om\big(A^{z}_{\a_{z}}
\big(A^{z}_{\a_{z}}\big)^{*}
\big)}}\bigg\}_{\a_{z}\in I_{z}}
$$ 
is an orthonormal basis for $E_{\bar z}\ch$ whenever 
$\{A^{z}_{\a_{z}}\Om\}_{\a_{z}\in I_{z}}$ is an orthonormal basis for 
$E_{z}\ch$.
\end{prop}
\begin{proof}
The fact that $M_z\Om$ is dense in $E_{z}\ch$ follows from 
Proposition 3.2 of \cite{NSZ}. Then, by exchanging the role 
between $M$ and $M'$, $(M')_z\Om$ is also dense in $E_{z}\ch$. By 
taking into account that the adjoint action 
of $U$ on $M$ (or equivalently on $M'$) is an automorphism, if $A\in 
M_z$ is nonnull, $A^{*}$ is a nonnull element of $M_{\bar z}$, that 
is $z\in\s_{\mathop{\rm pp}}(U)$ implies $\bar z\in\s_{\mathop{\rm pp}}(U)$ 

For $z\in\s_{\mathop{\rm pp}}(U)$, choose an 
orthonormal basis
$\{u^{z}_{\a_{z}}\}_{\a_{z}\in I_{z}}\subset M_z\Om$ 
for $M_z\Om$ which always exists by Zorn Lemma (see 
e.g. \cite{TL}, Section II.6). By the Parseval identity, 
a generic element $x$ in the completion of $M_z\Om$, the last 
coinciding with $E_{z}\ch$, is written as 
$$
x=\sum_{\a_{z}\in I_{z}}a_{\a_{z}}u^{z}_{\a_{z}}\,,
$$
where $\{a_{\a_{z}}\}_{\a_{z}\in I_{z}}\subset\bc$ is any  
square--summable sequence. 

The last property follows as $S\lceil_{M_z\Om}$ 
($S^{*}\lceil_{(M')_z\Om}$) realizes an algebraic antilinear isomorphism between
$M_z\Om$ and $M_{\bar z}\Om$ ($(M')_z\Om$ and $(M')_{\bar z}\Om$), 
$S$ being the Tomita involution of $M$ associated to the standard 
vector $\Om$ (cf. \cite{SZ}).
\end{proof}
\begin{prop}
\label{gnsc1}
Let the ergodic $C^*$--dynamical system $\big(\ga,\a,\om\big)$ be such that 
$s(\om)$ is central. 

Then, with the previous notations, $E_{z}\ch$ is 
one dimensional for each $z\in\s_{\mathop{\rm pp}}(U)$, and we can 
choose its single generator as $V_{z}\Om$ ($W_{z}\Om$) for some 
unitary $V_{z}\in M$ ($W_{z}\in M'$). Finally, $z\in\s_{\mathop{\rm pp}}(U)$ is 
a subgroup of $\bt$.
\end{prop}
\begin{proof}
Let $z\in\s_{\mathop{\rm pp}}(U)$ and choose nonnull elements $A,B\in M_{z}$
which always exist by Proposition 3.2 of \cite{NSZ}. Then
$A^{*}B=\a I$ and $BA^{*}=\b I$  for some nonnull numbers $\a,\b$. 
Indeed, $A^{*}B\Om$ is invariant under $U$. Thus by ergodicity,
$A^{*}B\Om=\a\Om$, and by the fact that $\Om$ is separating, $A^{*}B=\a I$. 
In addition, suppose that $\a=0$. As $AA^{*}$ is a nonnull multiple, 
say $c$, of 
the identity, we have $AA^{*}B=0$, which means $B=0$, a 
contradiction. At the same way, we verify $BA^{*}\neq0$. 
Now, $\a^{-1}A^{*}$ and $\b^{-1}A^{*}$ are left and right inverses of 
$B$. This means that $B$ is invertible and $B^{-1}=\a^{-1}A^{*}$. At 
the same way, $A$ is invertible too. Moreover,
$AB^{-1}=\a^{-1}AA^{*}=\a^{-1}cI$. This means $A=\a^{-1}cB$, that is 
$A$ is a multiple of $B$. In 
addition, in this situation $AA^{*}=cI$ means that $A$ is a multiple 
of the unitary $A/\sqrt{c}$. If $z,w\in\s_{\mathop{\rm pp}}(U)$, let 
$V_{z}\in M_{z}$, $V_{w}\in M_{w}$ be the corresponding unitaries. 

Then $V_{z}V_{w}\in M_{zw}$ is nonnull, that is $zw\in\s_{\mathop{\rm pp}}(U)$.
\end{proof}

\section{the entangled ergodic theorem in the almost periodic case}

The present section is devoted to the almost 
periodic situation, without any restriction relative to the operators 
appearing in \eqref{0}, and the pair--partition 
$\a:\{1,\dots,2k\}\mapsto\{1,\dots,k\}$. In this way, we improve the results 
in Section 3 of \cite{F} where only very special pair--partitions
were considered. 
We start by recalling for the reader convenience the known results 
relative to the entangled ergodic theorem. 

Let $U$ be 
a unitary operator acting on the Hilbert space $\ch$, and for $m\geq k$, 
$\a:\{1,\dots,m\}\mapsto\{1,\dots, k\}$ a partition of the 
set $\{1,\dots,m\}$.\footnote{A partition $\a:\{1,\dots,m\}\mapsto\{1,\dots, k\}$
of the set made of $m$ elements in $k$ parts is nothing but a surjective 
map, the parts of $\{1,\dots,m\}$ being the preimages 
$\{\a^{-1}(\{j\})\}_{j=1}^{k}$.} It was shown in Theorem 2.6 of \cite{F} 
that the multiple Cesaro mean in \eqref{0} converges in the weak 
operator topology when $A_{1},\dots,A_{m-1}\in\ck(\ch)$, $\ck(\ch)$ 
being the algebra of all the compact operators acting on $\ch$. In the case 
of a pair--partition $\a:\{1,\dots,2k\}\mapsto\{1,\dots, k\}$, we 
have (cf. \cite{F}, Theorem 2.5),
\begin{align*}
\mathop{\rm w\!-\!lim}_{N\to+\infty}&
\bigg\{\frac{1}{N^{k}}\sum_{n_{1},\dots,n_{k}=0}^{N-1}
U^{n_{\a(1)}}A_{1}U^{n_{\a(2)}}\cdots 
U^{n_{\a(2k-1)}}A_{2k-1}U^{n_{\a(2k)}}\bigg\}\\
=&S_{\a;A_{1},\dots,A_{2k-1}}\,,
\end{align*} 
where $A_{1},\dots,A_{m-1}\in\ck(\ch)$, and
$S_{\a;A_{1},\dots,A_{2k-1}}$ 
is given in \eqref{symb}. After passing to the finite rank operators, 
the Cesaro Mean in \eqref{0} disentangles, and the proof 
follows by Lebesgue dominated convergence theorem.

In order to treat the almost periodic case, from now on we suppose in the 
present section  
that $\ch$ is generated by 
the eigenvectors of $U$. 

The proof of the following result relies upon 
the mean ergodic theorem \eqref{jvn}, by 
showing that one can reduce oneself to the dense 
subspace algebraically generated by the eigenvectors of $U$.
\begin{thm}
\label{qper2}    
Let $U$ be an almost periodic unitary 
operator acting on the Hilbert space $\ch$. Then 
for each pair--partition $\a:\{1,\dots,2k\}\mapsto\{1,\dots,k\}$,
and $A_{1},\dots,A_{2k-1}\in\cb(\ch)$,
\begin{align}
\label{sta1}
\mathop{\rm s\!-\!lim}_{N\to+\infty}&\bigg\{\frac{1}{N^{k}}\sum_{n_{1},\dots,n_{k}=0}^{N-1}
U^{n_{\a(1)}}A_{1}U^{n_{\a(2)}}\cdots 
U^{n_{\a(2k-1)}}A_{2k-1}U^{n_{\a(2k)}}\bigg\}\nn\\
=&S_{\a;A_{1},\dots,A_{2k-1}}\,.
\end{align} 
\end{thm}
\begin{proof}
We treat the case of the partition $\{1,2,1,3,2,3\}$, the general 
case follows analogously.
Fix $\eps>0$, and suppose that $A,B,C,D,F\in\cb(\ch)$ have norm  
one. Let $I_{\eps}$ be such that 
$$
\bigg\|x-\sum_{\s\in I_{\eps}}E_{\s}x\bigg\|<\eps\,.
$$
For each $\s\in I_{\eps}$, let $I_{\eps}(\s)$ be such that 
$$
\bigg\|FE_{\s}x-\sum_{\t\in 
I_{\eps}(\s)}E_{\t}FE_{\s}x\bigg\|<\frac{\eps}{|I_{\eps}|}\,.
$$
By taking into account the mean ergodic theorem \eqref{jvn}, choose $N_{\eps}$ such that 
$$
\bigg\|\bigg(\frac{1}{N}\sum_{n=0}^{N-1}(\s U)^{n}-E_{\bar\s}\bigg)
DE_{\t}FE_{\s}x\bigg\|<\frac{\eps}{{\displaystyle\sum_{\s\in I_{\eps}}|I_{\eps}(\s)|}}\,,
$$
whenever $N>N_{\eps}$ and $\s\in I_{\eps}$, $\t\in I_{\eps}(\s)$.  
Finally, for each $\s\in I_{\eps}$, $\t\in I_{\eps}(\s)$, let $I_{\eps}(\s,\t)$
be such that 
$$
\bigg\|CE_{\bar\s}DE_{\t}FE_{\s}x-\sum_{\r\in I_{\eps}(\s,\t)}E_{\r}C
E_{\bar\s}DE_{\t}FE_{\s}x\bigg\|<\frac{\eps}
{{\displaystyle\sum_{\s\in I_{\eps}}|I_{\eps}(\s)|}}\,.
$$
We obtain by \eqref{zaz},
\begin{align*}
&\bigg\|\frac{1}{N^{3}}\sum_{k,m,n=0}^{N-1}U^{k}AU^{m}BU^{k}CU^{n}DU^{m}FU^{n}x
-S_{\a;A,B,C,D,F}x\bigg\|\\
\leq&5\eps
+\sum_{\s\in I_{\eps}}\sum_{\t\in I_{\eps}(\s)}\sum_{\r\in I_{\eps}(\s,\t)}
\bigg\|\bigg(\frac{1}{N}\sum_{k=0}^{N-1}(\r U)^{k}\bigg)A
\bigg(\frac{1}{N}\sum_{m=0}^{N-1}(\t U)^{m}\bigg)\\
\times&BE_{\r}C
E_{\bar\s}DE_{\t}FE_{\s}x-E_{\bar\r}AE_{\bar\t}B
E_{\r}CE_{\bar\s}DE_{\t}FE_{\s}x\bigg\|\,.
\end{align*}
Taking the limsup on both sides, we obtain the assertion by the mean ergodic theorem 
\eqref{jvn}, 
by taking into account the fact that $\eps>0$ is arbitrary.
\end{proof}

\section{an ergodic theorem for non invariant states}

The present section concerns the generalization to the quantum 
case of a well--known classical ergodic theorem due to H. Furstenberg,
for non invariant measures (cf. \cite{Fu, Fu1}). Such a theorem has a 
natural application to the noncommutative case of ``diagonal measures''.

We start with a $C^*$--dynamical system $\big(\ga,\a,\f\big)$, 
whose GNS covariant representation 
is denoted as $\big(\ch_{\f},\pi_{\f},U,\F\big)$. 
Consider another state
$\om\in\cs(\ga)$. Notice that $\om$ is supposed in general to be neither invariant 
under the action of $\a$, nor normal w.r.t. $\f$. Let 
$\big(\ch_{\om},\pi_{\om},\Om\big)$ be the GNS representation of $\om$. 
Let $\gb$ be a 
$*$--subalgebra of $\ga$.  

The following definition is nothing but the 
generalization of Definition 4.4 of \cite{Fu1} to the quantum case. 
\begin{defin}
\label{fu}
The state $\om$ is said to be generic for $\big(\ga,\a,\f\big)$ 
w.r.t. $\gb$ if for each $B\in\gb$,
$$
\lim_{N\to+\infty}\frac{1}{N}\sum_{n=0}^{N-1}\om(\a^{n}(B))=\f(B)\,.
$$
\end{defin}

Let $E_{1}$ be the selfadjoint projection onto the invariant vectors 
for the unitary $U$. Suppose that $\pi_{\f}(\gb)\F\bigcap E_{1}\ch_{\f}$ 
is dense in $E_{1}\ch_{\f}$.
\begin{lem}
\label{fu2}
Under the above conditions, 
\begin{equation}
\label{fu4}
\pi_{\f}(B)\F\in\pi_{\f}(\gb)\F\bigcap E_{1}\ch_{\f}\mapsto \pi_{\om}(B)\Om\in\ch_{\om}
\end{equation}
uniquely defines a partial isometry $V:\ch_{\f}\mapsto\ch_{\om}$ such 
that $V^{*}V=E_{1}\ch_{\f}$. 
\end{lem}
\begin{proof}
It is enough to show that, under our assumptions, the map in \eqref{fu4} is isometric.
We get for each $B\in\gb$ invariant under 
the action of $\a$, first by taking into account the invariance of 
$B$, and then the genericity of $\om$,
\begin{align*}
&\|\pi_{\om}(B)\Om\|^{2}\equiv\om(B^{*}B)
\equiv\frac{1}{N}\sum_{n=0}^{N-1}\om(\a^{n}(B^{*}B))\\
\equiv&\lim_{N\to+\infty}\frac{1}{N}\sum_{n=0}^{N-1}\om(\a^{n}(B^{*}B))
=\f(B^{*}B)\equiv\|\pi_{\f}(B)\F\|^{2}\,.
\end{align*}
\end{proof}    
For $A\in\cb(\ch)$, denote $|A|:=(A^{*}A)^{1/2}$.
We need also the following technical results.
\begin{lem}
\label{tec}
Let $A_{1},\dots,A_{n}$ be bounded operators acting on the Hilbert space 
$\ch$. Then
$$
\bigg|\frac{1}{n}\sum_{k=1}^{n}A_{k}\bigg|^{2}\leq
\frac{1}{n}\sum_{k=1}^{n}|A_{k}|^{2}\,.
$$
\end{lem}
\begin{proof}
The proof easily follows if one verifies
$A^{*}B+B^{*}A\leq A^{*}A+B^{*}B$. But,
$$
0\leq(A-B)^{*}(A-B)= A^{*}A+B^{*}B-(A^{*}B+B^{*}A)\,.
$$
\end{proof}
\begin{lem}
\label{tec1}
Let $\{a_{k}\}_{k\in\bn}\subset\bc$ be a bounded sequence. 
Then for each fixed $M$, 
$$
\lim_{N\to+\infty}\bigg(\frac{1}{MN}\sum_{n=0}^{N-1}\sum_{m=0}^{M-1}a_{m+n}
-\frac{1}{N}\sum_{n=0}^{N-1}a_{n}\bigg)
$$
at a rate depending on $\sup_{k}|a_{k}|$.
\end{lem}
\begin{proof}
The proof follows by taking into account    
$$
\bigg|\frac{1}{MN}\sum_{n=0}^{N-1}\sum_{m=0}^{M-1}a_{m+n}
-\frac{1}{N}\sum_{n=0}^{N-1}a_{n}\bigg|\leq\frac{(M-1)(M+2)}{MN}\sup_{k}|a_{k}|\,.
$$
\end{proof}
The following theorem is nothing but the announced
generalization of Theorem 4.14 of \cite{Fu1} (see also \cite{Fu}) 
to the quantum case.
\begin{thm}
\label{fu1}
Let $\om$ be generic for $\big(\ga,\a,\f\big)$ 
w.r.t. a $*$--subalgebra $\gb$ which is globally stable under the 
action of $\a$, and satisfies 
$$
\overline{\pi_{\f}(\gb)\F\bigcap E_{1}\ch_{\f}}=E_{1}\ch_{\f}\,.
$$
Then for each $B\in\gb$,
$$
\lim_{N\to+\infty}\frac{1}{N}\sum_{n=0}^{N-1}\pi_{\om}(\a^{n}(B))\Om
=V\left(\pi_{\f}(B)\F\right)\,,
$$
$V:\ch_{\f}\mapsto\ch_{\om}$ being the partial isometry given in 
Lemma \ref{fu2}.
\end{thm}
\begin{proof}
Let $B\in\gb$ and $\eps>0$ be given. Choose a 
$\a$--invariant $B_{\eps}\in\gb$ such that 
$\|\big(E_{1}\pi_{\f}(B)-\pi_{\f}(B_{\eps})\big)\F\|\leq\eps$. By the 
mean ergodic theorem \eqref{jvn},
$$
\frac{1}{M}\sum_{m=0}^{M-1}\pi_{\f}(\a^{m}(B-B_{\eps}))\big)\F
\longrightarrow\big(E_{1}\pi_{\f}(B)-\pi_{\f}(B_{\eps})\big)\F\,.
$$
Thus, for $M$ sufficiently large,
$$
\f\bigg(\bigg|\frac{1}{M}\sum_{m=0}^{M-1}\a^{m}(B-B_{\eps})\bigg|^{2}\bigg)<\eps^{2}\,.
$$
Denote 
$$
\G:=\frac{1}{M}\sum_{m=0}^{M-1}\a^{m}(B-B_{\eps})\,.
$$
By hypotesis, $\G^{*}\G\in\gb$, and as $\om$ is generic w.r.t. $\gb$,
$$
\frac{1}{N}\sum_{n=0}^{N-1}\om(\a^{n}(\G^{*}\G))
\longrightarrow\f(\G^{*}\G)\,.
$$
So, for each $N$ sufficiently large, 
$$
\frac{1}{N}\sum_{n=0}^{N-1}\om(\a^{n}(\G^{*}\G))<\eps^{2}\,.
$$
By applying Lemma \ref{tec}, we have for a fixed $M$ sufficiently 
large
and each $N$ sufficiently large,
\begin{equation}
\label{lm}
\om\bigg(\bigg|\frac{1}{NM}\sum_{n=0}^{N-1}\sum_{m=0}^{M-1}
\a^{n+m}(B-B_{\eps})\bigg|^{2}\bigg)<\eps^{2}\,.
\end{equation}
By taking into account Lemma \ref{tec1}, \eqref{lm} becomes
$$
\om\bigg(\bigg|\frac{1}{N}\sum_{n=0}^{N-1}
\a^{n}(B)-B_{\eps}\bigg|^{2}\bigg)<\eps^{2}
$$
which means
$$
\bigg\|\frac{1}{N}\sum_{n=0}^{N-1}
\pi_{\om}(\a^{n}(B)-B_{\eps})\Om\bigg\|<\eps
$$
for each large $N$. Thus, we obtain for each sufficiently large $N$,
\begin{align*}
&\bigg\|\frac{1}{N}\sum_{n=0}^{N-1}
\pi_{\om}(\a^{n}(B))\Om-VE_{1}\pi_{\f}(B)\F\bigg\|
\leq\bigg\|\frac{1}{N}\sum_{n=0}^{N-1}
\pi_{\om}(\a^{n}(B)-B_{\eps})\Om\bigg\|\\
+&\big\|V\big(\pi_{\f}(B_{\eps})-E_{1}\pi_{\f}(B)\big)\F\big\|
<\eps+\big\|\big(\pi_{\f}(B_{\eps})-E_{1}\pi_{\f}(B)\big)\F\big\|<2\eps\,.
\end{align*}
\end{proof}

\section{the case of ``diagonal measures''}

The present section is devoted to the natural generalization to the 
quantum case of the celebrated result due to H. Furstenberg relative to 
the diagonal measures (cf. \cite{Fu}, see also \cite{Fu1}, Section 
4.4).

We start with a $C^*$--dynamical system $\big(\ga,\a,\om\big)$, 
together with its GNS covariant representation $\big(\ch,\pi,U,\Om\big)$. 
Denote $M:=\pi(\ga)''$,  
the von Neumann algebra acting on $\ch$ generated by 
the representation $\pi$. The commutant von Neumann algebra is denoted as $M'$.
Suppose further that the support $s(\om)$ in $\ga^{**}$
is central. 

Let $\gam:=M\ots{max}M'$ be the completion of the algebraic 
tensor product $\gn:=M\otimes M'$ w.r.t. the maximal $C^*$--norm (cf. 
\cite{T}, Section IV.4). It is easily seen that on $\gam$ the 
following two states are automatically well--defined. The first one is 
the canonical product state 
$$
\f(A\otimes B):=\langle A\Om,\Om\rangle\langle B\Om,\Om\rangle\,,
\quad A\in M\,,\,\,B\in M'\,.
$$
The second one is uniquely defined by
$$
\psi(A\otimes B):=\langle AB\Om,\Om\rangle\,,
\quad A\in M\,,\,\,B\in M'\,.
$$

The state $\psi$ can be considered the (quantum analogue of the) ``diagonal 
measure'' of the ``measure'' $\f$.

On $\gam$ is also uniquely defined the automorphism
$$
\g:=\ad{}\!_{U}\otimes\ad{}\!_{U^{2}}\,,
$$
see \cite{T}, Proposition IV.4.7. Of course, $\big(\gam,\g,\f\big)$ is 
a $C^*$--dynamical system whose GNS covariant representation is 
precisely 
$\big(\ch\oots\ch,\id\otimes\id,U\otimes U^{2},\Om\otimes\Om\big)$. 
Denote $E_{1}$ the selfadjoint 
projection onto the invariant vectors under $U\otimes U^{2}$.
Notice that the 
$*$--subalgebra $\gn$ is globally stable under the action of $\g$.

In addition, again by Proposition IV.4.7 of \cite{T}, 
$$
\s(A\otimes B):=AB\,,
\quad A\in M\,,\,\,B\in M'\,.
$$
uniquely defines a representation of $\gam$ on $\ch$ such that 
$\big(\ch,\s,\Om\big)$ is precisely the GNS representation of the 
state $\psi$.\footnote{Notice that, even if it is enough for our 
pourpose to consider $M\ots{max}M'$,
all these properties hold true 
for $M\ots{bin}M'$, the latter being the 
completion of $M\otimes M'$ with the binormal $C^*$--norm (cf. 
\cite{EL}).}
\begin{prop}
\label{diam}
Suppose that $\big(\ga,\a,\om\big)$ is ergodic. Then the state $\psi\in\cs(\gam)$ is 
generic for $\big(\gam,\g,\f\big)$ w.r.t. $\gn$.
\end{prop}
\begin{proof}
Let $A\in M$, $B\in M'$. Then by the mean ergodic theorem \eqref{jvn},
\begin{align*}
\frac{1}{N}\sum_{n=0}^{N-1}\psi(\g^{n}(A&\otimes B))
=\frac{1}{N}\sum_{n=0}^{N-1}\langle AU^{n}B\Om,\Om\rangle\\
&\equiv\bigg\langle A\bigg(\frac{1}{N}\sum_{n=0}^{N-1}U^{n}\bigg)B\Om,
\Om\bigg\rangle\\
&\longrightarrow\langle A\Om,\Om\rangle\langle B\Om,\Om\rangle
\equiv\f(A\otimes B)\,.
\end{align*}
\end{proof}
\begin{thm}
\label{diam1}    
Let $\big(\ga,\a,\om\big)$ be an ergodic $C^*$--dynamical system such 
that its support $c(\om)$ in $\ga^{**}$ is central. Then 
with the above notations, the following assertions hold true.
\begin{itemize}
\item[(i)] Let ${\displaystyle\sum_{j}A_{j}\otimes B_{j}}$ be the 
generic element of $\gn$. The map 
$$
\sum_{j}A_{j}\Om\otimes B_{j}\Om\in\gn\Om\bigcap E_{1}\ch\oots\ch
\mapsto\sum_{j}A_{j}B_{j}\Om\in\ch
$$
uniquely defines a partial isometry $V:\ch\oots\ch\mapsto\ch$ such 
that $V^{*}V=E_{1}\ch\oots\ch$.
\item[(ii)] For $A\in M$, $B\in M'$,
$$
\lim_{N\to+\infty}\frac{1}{N}\sum_{n=0}^{N-1}U^{n}AU^{n}B\Om
=V\left(A\Om\otimes B\Om\right)\,.
$$
\end{itemize}
\end{thm}
\begin{proof}
Define ${\displaystyle\S:=\{(z,w)\in\s_{\mathop{\rm pp}}
(U)\times\s_{\mathop{\rm pp}}(U)\,:\,
zw^{2}=1\}}$.

Then by Lemma 4.18 of \cite{Fu1},
$$
E_{1}=\bigoplus_{s\in\S}E^{U}_{z_{s}}\otimes E^{U}_{w_{s}}\,,
$$
$E^{U}_{z}$ being the selfadjoint projection onto 
the eigenspace of $U$ corresponding to the eigenvalue $z$. As $U$ is 
ergodic, by Proposition \ref{gnsc1} $E^{U}_{z}\ch$ is one dimensional, 
and $E^{U}_{z_{s}}\ch$ and $E^{U}_{w_{s}}\ch$ are generated by 
$V_{z_{s}}\Om$, $W_{w_{s}}\Om$, where $V_{z_{s}}$ and $W_{w_{s}}$ are 
unitaries of $M_{z_{s}}$, $(M')_{w_{s}}$ respectively. Thus,
$E^{U}_{z_{s}}\ch\oots E^{U}_{w_{s}}\ch$ is one dimensional, and it 
is generated by $V_{z_{s}}\Om\otimes W_{w_{s}}\Om$. This means that 
$\gn\Om\bigcap E_{1}\ch\oots\ch$ is dense in $E_{1}\ch\oots\ch$. 

The 
assertions follow from Proposition \ref{diam} and
Theorem \ref{fu1}, by taking into account that 
$\gn$ is left globally 
invariant by $\ad{}\!_{U}\otimes\ad{}\!_{U^{2}}$.
\end{proof}

The following results are a direct consequence of the 
previous one.
\begin{cor}
\label{cccc}
Under the hypotheses of Theorem \ref{diam1},
\begin{equation}
\label{sta}
\mathop{\rm s\!-\!lim}_{N\to+\infty}
\frac{1}{N}\sum_{n=0}^{N-1}U^{n}AU^{n}
=V\left(A\Om\otimes\,\cdot\,\right)\,,
\end{equation}
for each $A\in M\bigcup M'$.
\end{cor}
\begin{proof}
By a standard density argument, if $A\in M$, the proof follows from Theorem 
\ref{diam1}, by taking into account that $\Om$ is cyclic from $M'$. 
The proof for $A\in M'$ follows by exchanging the role between $M$ 
and $M'$. 
\end{proof}

By taking into account Proposition \ref{gnsc} and Proposition 
\ref{gnsc1}, it is straightforward to verify that \eqref{sta} 
for $A\in M\bigcup M'$ coincides with \eqref{sta1} when $U$ is almost 
periodic and ergodic.\footnote{Notice that in this situation,
$\s_{\mathop{\rm pp}}(U)=\s_{\mathop{\rm pp}}^{\mathop{\rm a}}(U)$,
with $\s_{\mathop{\rm pp}}^{\mathop{\rm a}}(U)$ given in \eqref{asy}, 
and $\s(U)=\overline{\s_{\mathop{\rm pp}}(U)}$.}
Namely, 
$$
E_{z}=\langle\,\cdot\,,V_{z}\Om\rangle
V_{z}\Om
=\langle\,\cdot\,,W_{z}\Om\rangle W_{z}\Om
$$
for unitaries $V_{z}$, $W_{z}$ in 
$M_{z}$, $(M')_{z}$ respectively. Then
\begin{align*}
=&\sum_{w\in\s_{\mathop{\rm pp}}(U)}E_{\bar w}AE_{w}\xi
=\sum_{w\in\s_{\mathop{\rm pp}}(U)}\sum_{\{z\in\s_{\mathop{\rm pp}}
(U)\,:\,zw^{2}=1\}}
E_{zw}AE_{w}\xi\\
\equiv&\sum_{\{z,w\in\s_{\mathop{\rm pp}}(U)\,:\,zw^{2}=1\}}E_{zw}AE_{w}\xi\\
=&\sum_{\{z,w\in\s_{\mathop{\rm pp}}(U)\,:\,zw^{2}=1\}}\langle\xi,W_{w}\Om\rangle
\langle AW_{w}\Om,V_{z}W_{w}\Om\rangle V_{z}W_{w}\Om\\
=&\sum_{\{z,w\in\s_{\mathop{\rm pp}}(U)\,:\,zw^{2}=1\}}\langle A\Om,V_{z}\Om\rangle
\langle\xi,W_{w}\Om\rangle
V_{z}W_{w}\Om\\
=&V\left(A\Om\otimes\xi\right)\,.
\end{align*}    
\begin{cor}
Under the hypotheses of Theorem \ref{diam1},
\begin{align*}
&\lim_{N\to+\infty}\frac{1}{N}\sum_{n=0}^{N-1}\om\left(A_{0}\a^{n}(A_{1})
\a^{2n}(A_{2})\right)\\
&=\left\langle V\left(\pi(A_{1})\Om\otimes\pi(A_{2})
\Om\right),\pi(A^{*}_{0})\Om\right\rangle\,.
\end{align*}
\end{cor}
\begin{proof}
A simple application of Corollary \ref{cccc}.
\end{proof}
Suggested by the Abelian situation (cf \cite{Fu1}, pag. 96), one can ask 
for the convergence of the Cesaro mean
\begin{equation}
\label{cm}
\frac{1}{N}\sum_{n=0}^{N-1}\pi\big(\a^{nm_{1}}(A_{1})
\a^{nm_{2}}(A_{2})\big)\Om   
\end{equation}
for the other cases with fixed $0<m_{1}<m_{2}$. Starting from 
\begin{align}
\label{dma}
&\lim_{N\to+\infty}\frac{1}{N}\sum_{n=0}^{N-1}\left\langle AU^{n(m_{2}-m_{1})}B
\Om,\Om\right\rangle\nn\\
&=\sum_{\{z\in\s_{\mathop{\rm pp}}(U)\,:\,z^{m_{2}-m_{1}}=1\}}
\left\langle AE_{z}B\Om,\Om\right\rangle
\end{align}
whenever $A\in\pi(\ga)''$, $B\in\pi(\ga)'$, in order to apply Theorem
\ref{fu1} one firstly demand if \eqref{dma} uniquely defines a 
state on $\pi(\ga)''\ots{max}\pi(\ga)'$. Such a state will be
necessarily invariant under 
the action of $\ad{}\!_{U^{m_{1}}}\otimes\ad{}\!_{U^{m_{2}}}$. This is certainly 
true when $\#\{z\in\s_{\mathop{\rm pp}}(U)\,:\,z^{m_{2}-m_{1}}=1\}$ is finite.

After verifying the remaining hypotheses of Theorem
\ref{fu1}, one might argue that the Cesaro mean in \eqref{cm} 
converges, at least when the subspace consisting of all the invariant 
vectors for $U^{m_{2}-m_{1}}$ is finite dimensional.
We end with the simplest case of weakly mixing dynamical systems. Then 
we have an alternative proof of (a weaker result than) Theorem 1.3 of
\cite{NSZ}, following the line of Theorem \ref{diam1}.
\begin{prop}
\label{diamm1}    
Let $\big(\ga,\a,\om\big)$ be a weakly mixing $C^*$--dynamical system such 
that its support $c(\om)$ in $\ga^{**}$ is central, and 
$0<m_{1}<m_{2}$ natural numbers. 

Then with the above 
notations,
\begin{equation*}
\mathop{\rm s\!-\!lim}_{N\to+\infty}
\frac{1}{N}\sum_{n=0}^{N-1}U^{nm_{1}}AU^{nm_{2}}
=\langle A\Om,\Om\rangle\langle\,\cdot\,,\Om\rangle\Om\,,
\end{equation*}
for each $A\in M\bigcup M'$.
\end{prop}
\begin{proof}
We apply Theorem \ref{fu1} by considering
$$
\g:=\ad{}\!_{U^{m_{1}}}\otimes\ad{}\!_{U^{m_{2}}}\,.
$$

Indeed, $E_{1}\ch\oots\ch=\bc\Om\otimes\Om$ where
$E_{1}$ is the selfadjoint 
projection onto the invariant vectors under $U^{m_{1}}\otimes 
U^{m_{2}}$. In addition, if $A\in M$, $B\in M'$,
\begin{align*}
\frac{1}{N}\sum_{n=0}^{N-1}\psi(\g^{n}(A&\otimes B))
=\frac{1}{N}\sum_{n=0}^{N-1}\langle AU^{n(m_{2}-m_{1})}B\Om,\Om\rangle\\
&\longrightarrow\langle A\Om,\Om\rangle\langle B\Om,\Om\rangle
\equiv\f(A\otimes B)\,.
\end{align*}
\end{proof}

\end{document}